\documentclass{amsbook}
\usepackage{amssymb}
\usepackage{amsfonts}

\makeatletter
\numberwithin{equation}{section}
\numberwithin{figure}{section}
\numberwithin{equation}{section}
\numberwithin{figure}{section}
\numberwithin{equation}{section}

\theoremstyle{definition}

\newtheorem*{acknowledgements*}{Acknowledgements}
\theoremstyle{remark}

\numberwithin{theorem}{section}   
\input{tcilatex}

\begin{document}
\frontmatter
\title{$H^{\infty }$ and Complex Interpolation.}
\author{Jaak Peetre}
\address{Lund Institute of Technology, Lund, Sweden}
\author{Per G. Nilsson (Typist)}
\address{Nilsson: Stockholm, Sweden}
\email{pgn@plntx.com}
\date{\today }

\begin{abstract}
This note is an (exact) copy of the report of Jaak Peetre, "$H^{\infty }$
and Complex Interpolation". Published as Technical Report, Lund (1981). Some
more recent general references have been added, some references updated
though (in \textit{italics}) and some misprints corrected.
\end{abstract}

\maketitle

\begin{center}
\bigskip
\end{center}

\section*{${\protect\Huge H}^{\infty }$ and complex interpolation.}

\section*{Jaak Peetre}

\section*{0. Introduction}

\pagenumbering{arabic}

In the theory of complex interpolation one usually defines the interpolation
spaces in question using holomorphic functions taking continuous boundary
values in the classical sense. (For the complex interpolation the main
source is still Calderon's classic, albeit (\underline{helas }!) tough
reading paper \cite{5-Cal}. For an excellent introduction we recommend chap.
4 in the book \cite{3-BeL0}.)\footnote{%
For references, post 1981, see notably Brudnyi-Krugljak \cite[Chap. 2, Chap.
4]{27-BK}, Brudnyi-Krein-Semenov \cite{28-BrKrSe}, Ovchinnikov \cite{29-Ov84}%
, Lindemulder-Lorist \cite{30-LiFl24} and the references listed.} That is,
we have to deal with the space $A$.$\medskip $

The main thesis advanced in this paper is that in many questions is much
more advantageous to consider the space $H^{\infty }$ (of bounded
holomorphic functions) and further that it is in this connection often
easier to work with distributional boundary values rather than with
classical pointwise (a.e.) ones. $\medskip $

Our main result in this direction (sec. 4, theorem) says that one gets the
same interpolation space $\overline{X}_{\theta }$ no matter whether one uses 
$A$ or $H^{\infty }$. (Here and in the sequel $\overline{X}=\left(
X_{0},X_{1}\right) $ denotes any fixed (compatible) pair of Banach spaces;
see \cite{3-BeL0}.) We likewise obtain (sec. 5, theorem) a similar
characterization of the second Calderon space $\overline{X}^{\theta }$.$%
\medskip $

Our new description of $\overline{X}^{\theta }$ is formally obtained from
the one for $\overline{X}_{\theta }$ by substituting the functor $\hom
\left( L^{1},-\right) $ for the functor $L^{\infty }\left( -\right) $. As a
consequence we obtain the result (sec. 6, proposition) that if either $X_{0}$
or $X_{1}$ has the RN (Radon-Nikodym) property (see e.g. \cite{11-DiUh})
then $\overline{X}_{\theta }=\overline{X}^{\theta }$. Calderon \cite{5-Cal}
proved this with "reflexive" instead of "RN".$\medskip $

In fact our new proof of his result is almost trivial, modulo that fact that
"reflexive" entails "RN". We likewise give a new proof of the Calderon's
duality theorem \cite{5-Cal} and of Janson's recent characterization \cite%
{23-Jan} of the complex interpolation spaces in terms of orbit spaces. $%
\medskip $

Apart from these applications to essentially known results of the Calderon
theory we can use our new insight to treat \underline{vector valued}
versions of various problems traditionally considered for scalar $H^{\infty }
$. In particular we consider a natural analogue of the (free) interpolation
problem: Given a sequence of points $\left\{ z_{n}\right\} _{n=1}^{\infty }$
in the strip $\left\{ 0<Rez<1\right\} $ of the complex plane, which is
uniformly separating in the sense of Carleson \cite{5-Cal}. To find for a
given sequence $\left\{ w_{n}\right\} _{n=1}^{\infty }$ of elements of $%
X_{0}+X_{1}$, satisfying the appropriate boundedness condition, a function $%
f\in H^{\infty }\left( \overline{X}\right) $ (definition in sec. 3) such
that $f_{n}\left( z_{n}\right) =w_{n}$ ($n=1,2,...)$. Again this problem is
but a special case of a corresponding problem for the $\overline{\delta }$
operator (cf. \cite{21-Hor}). In this connection a recent paper by Jones 
\cite{24-Jon} was decisive. $\medskip $

Finally, to conclude this Introduction, we would like to point out what we
are really interested in is interpolation of \underline{infinitely many}
Banach spaces, not just two. (Concerning interpolation of infinitely many
spaces see the works of Coifman - Cwikel - Rochberg - Sagher - Weiss (the
"western" branch") and Krein - Nicolova (the "eastern" branch) respectively,
for instance the references \cite{8-CCRSW}, \cite{9-CCRSW}, \cite{25-KrNi}
listed here). It is mainly in order to keep the difficulties apart that we
presently have restricted ourselves to just \underline{two} spaces. We hope
however to be able to return to the general case in the near future. Let us
only mention here that the leading underlying idea, which has motivated us
in the present work too but which will find its true significance only in
these later developments, is that the true setting for the theory of
interpolation of infinitely many spaces is provided by \underline{Banach
bundles} in the sense of Fell \cite{15-Fel}, \cite{16-Fel}. (Thus what we
are really concerned with here is the Banach bundle valued $\overline{\delta 
}$ problem, rather than the vector valued one.) As a further indication of
that we are now firmly on the right track let us point out that whereas
Calderon's definition of $\overline{X}^{\theta }$ depends on special
features of the underlying domain (viz. the strip $\left\{ 0<Rez<1\right\} $%
), the invariance with respect to (vertical) translations, and does not
generalize, our's works a priori for any plane domain with a sufficiently
smooth boundary; indeed we can formally even consider generalizations with
complex manifolds (with a distinguished boundary) in any number of
dimensions (cf. Favini \cite{14-Fav}, Fernandez \cite{17-Fer}). $\medskip $

Having assimilated the particulars of this Introduction the reader can now
get a rough idea of how this paper is organized by glancing at the titles of
the various sections. $\medskip $

\underline{Acknowledgement}. My thanks are due to the following persons for
helpful suggestions, in particular what concerns bibliographic references,
in connection with this research: Lars G\aa rding, Svante Janson, Mario
Milman, Per Nilsson, Jan-Erik Roos, the last mentioned also for his
hospitably during the recent meeting of the Swedish Mathematical Society in
Stockholm (May 25-26, 1981).

\section*{1. The space $H^{\infty }\left( \Sigma \right) $.}

Let us denote by $S$ the (closed) strip $\left\{ 0\leq Rez\leq 1\right\} $
in the complex plane. Put further $S^{i}=\left\{ 0<Rez<1\right\} $, $\delta
S=\delta _{0}S\cup \delta _{1}S=\left\{ Rez=0\right\} \cup \left\{
Rez=1\right\} $.$\medskip $

We usually set $z=x+iy$.$\medskip $

Let $\overline{X}=\left( X_{0},X_{1}\right) $ be any (compatible) pair of
Banach spaces and write $\Sigma =\Sigma \overline{X}=X_{0}+X_{1}$; $\Sigma $
is a Banach space (cf. \cite{3-BeL0}, chap. 2). We say $f\in H^{\infty
}\left( \Sigma \right) $ if $f$ is bounded holomorphic function defined in $%
S^{i}$ with values in $\Sigma $. Clearly $H^{\infty }\left( \Sigma \right) $
is a Banach space in the norm $\left\Vert f\right\Vert =\sup_{z\in
S^{i}}\left\Vert f\left( z\right) \right\Vert $.$\footnote{%
See note $\left\langle 1\right\rangle .$}$ We denote by $H^{\infty }$ the
corresponding space of scalar (= complex valued) functions. $\medskip $

$\underline{\text{Claim.}}$ Boundary values of functions in $H^{\infty
}\left( \Sigma \right) $ exists in the distribution sense.$\medskip $

To substantiate this claim we first make a long detour and consider in all
generality the problem of existence of distributional boundary values. Sec.
2 will be devoted to it.

\section*{2. Digression on the existence of distributional boundary values.}

First we perform a rotation by an angle of $90^{\circ }$.$\medskip $

We begin with the \underline{scalar} case. Here we follow the most elegant
treatment of H\"{o}rmander \cite{22-Hor}. (The result in itself goes back to
antiquity (Archimedes?)$\medskip $.

Let thus $f$ be any (scalar) function holomorphic in an open "$\frac{1}{2}$
neighborhood" of $0\,$, that is, a set of the form $\left\{ \left\vert
z\right\vert <r,y>0\right\} $ ($r>0$). Assume that $f\left( z\right)
=O\left( y^{-A}\right) $ uniformly in $x$ on compact sets, for some $A\geq 0$%
.\footnote{%
See note $\left\langle 2\right\rangle .$}. Then $f$ can be continued to a
distribution in the full neighborhood $\left\{ \left\vert z\right\vert
<r\right\} $, still denoted by $f=f\left( z\right) $, with the support
contained in the "closed" $\frac{1}{2}$-neighborhood $\left\{ \left\vert
z\right\vert <r,y\geq 0\right\} $, such that for a suitable distribution $%
g=g\left( x\right) $ (in one variable) holds%
\begin{equation}
\int g\left( x\right) \psi \left( x,0\right) dx=-2\int_{y>0}f\left(
x+iy\right) \frac{\delta \psi \left( x,y\right) }{\delta \overline{z}}dxdy 
\tag{$\left( 1\right) $}  \label{eq-1}
\end{equation}%
for any (compactly supported) testfunction $\psi \left( x,y\right) $. (The
integrals in \ref{eq-1} are of course interpreted in distribution sense; $%
\frac{\delta }{\delta \overline{z}}=\frac{1}{2}\left( \frac{\delta }{\delta x%
}+i\frac{\delta }{\delta y}\right) $ is the Cauchy-Riemann operator.)$%
\medskip $

To prove this we first assume that $f$ takes continuous boundary values in
the classical sense (pointwise limits). Then \ref{eq-1} holds trivially with 
\underline{bona fide} integrals (just integrate by parts). Given any
testfunction $\phi =\phi \left( x\right) $ (in one variable we can determine 
$\psi $ with $\psi \left( x,0\right) =\phi \left( x\right) $ such that $%
\frac{\delta \psi \left( x,y\right) }{\delta \overline{z}}=O\left(
y^{A}\right) $. E.g. we take $\psi \left( x,y\right) =\phi \left( x\right) +%
\frac{iy}{1!}\phi ^{^{\prime }}\left( x\right) +...+\frac{\left( iy\right)
^{k}}{k!}\phi ^{\left( k\right) }\left( x\right) $ for $y$ close to $0$,
with $k\geq A$. Then we obtain from \ref{eq-1} the estimate%
\begin{equation}
\left\vert \int g\left( x\right) \phi \left( x\right) dx\right\vert \leq
C\left\Vert \phi \right\Vert _{k+1}  \tag{$\left( 2\right) $}  \label{eq-2}
\end{equation}%
where $\left\Vert \phi \right\Vert _{k+1}$ stands for a suitable norm
involving the maxima of the partial derivatives of $\phi $ up to order $k+1$
and $C$ a constant depending on the bound in $f\left( z\right) =O\left(
y^{-A}\right) $ only. The general case (distributional boundary values) is
now easily handled using \ref{eq-2} by replacing $f\left( z\right) $ by the
functions $f\left( z+i\epsilon \right) $ $\left( \epsilon >0\right) $ and
then letting $\epsilon $ tend to $0$.$\medskip $

In particular we see thus that $g\left( x\right) $ arises as the
distributional limit $\left( \epsilon \rightarrow 0\right) $ of the
functions $f\left( x+i\epsilon \right) $. (We therefore allow ourselves in
the sequel to write $f\left( x\right) $ instead of $g\left( x\right) $,
whenever this is convenient.) Notice also that it follows that any relation
between holomorphic functions can by continuity be translated into a similar
relation between the corresponding boundary values. (E.g. if $%
f_{1}^{^{\prime }}=f_{2}$ and $g_{1}$ and $g_{2}$ are the boundary values of 
$f_{1}$ and $f_{2}$ respectively then $g_{1}^{^{\prime }}=g_{2}$, in the
distribution sense of course.) We will somewhat pretentiously refer to this
property as the \underline{principle of permanence}.$\medskip $

We give now some additional comments on the above construction.

\underline{Remark 1.} It is clear that the "uniform" condition $f\left(
z\right) =O\left( y^{-A}\right) $ can be replaced by an integral condition
of the type continuous$\int \int \left\vert f\left( z\right) \right\vert
y^{A}dxdy<\infty $, with integration over smaller $\frac{1}{2}$%
-neighborhoods.

\underline{Remark 2.} For later use (sec. 8) we record that the argument
given also applies to functions $f$ satisfying the inhomogeneous
Cauchy-Riemann equation $\frac{\delta f}{\delta \overline{z}}=\mu $, $\mu $
a measure with finite mass on smaller $\frac{1}{2}-$neighborhoods.

\underline{Remark 3.} It is further easy to see that the condition $f\left(
z\right) =O\left( y^{-A}\right) $ (or its integral analogue; see remark 1)
values of a holomorphic function $f$ exists in the sense of \ref{eq-1} then
we must have $f\left( z\right) =O\left( y^{-A}\right) $ for some $A\geq 0$
uniformly in $x$ on compact sets.

\underline{Remark 4.} The proceeding discussion can formally be extended to
any partial differential operator with $C^{\infty }$ coefficients,
preferably elliptic though, even in $n$ variables. Notice that to the
special test function $\psi \left( x,y\right) $ constructed then corresponds
the formal solution of the Cauchy problem truncated after sufficiently many
steps.$\medskip $

We are now ready for the \underline{vector valued} case.$\medskip $

First of all we nail down that by a distribution with values in $\Sigma $
(or any other Banach space for that matter) we mean a continuous linear map $%
f$ from $\mathcal{D}$ (the space of all testfunctions) into $\Sigma $. For
the value $f\left( \phi \right) $ of $f$ at the element $\phi \in \mathcal{D}
$ we again use the notation $\int f\left( x\right) \phi \left( x\right) dx$,
in the case of one dimension, and similarly in two dimensions.$\medskip $

Let thus $f$ be a holomorphic function in the same $\frac{1}{2}$-
neighborhood, however with values in $\Sigma $. Assume that $\left\Vert
f\left( z\right) \right\Vert _{\Sigma }=O\left( y^{-A}\right) $ uniformly in 
$x$ on compacts sets, for some $A\geq 0$. Then the proceeding considerations
carry over \underline{mutatis mutandi}. We conclude that $f$ admits vector
values distributional boundary values $g$ in the sense that the inequality %
\ref{eq-1} holds true for any test function $\psi $.$\medskip $

Finally we state a simple lemma which will be needed in the sequel.

\underline{Lemma,} Assume that $g$ in \ref{eq-1} happens to be a continuous
function. Then $f$ has a continuous extension "up to the boundary" (i.e. to
the set $\left\{ z:\left\vert z\right\vert <r,y\geq 0\right\} $ and
consequently the boundary values exists in the classical sense.

\underline{Proof:} The proof is equally simple in the scalar as in the
vector valued case. It suffices to establish the continuity at the origin $0$%
. Introduce to this end a suitable cut-off function $\kappa $ equal to $1$
near $0$. Then $\kappa f$ clearly can be represented as the sum of a
function which is $C^{\infty }$ near $0$ and the Poisson integral of $\kappa
g$. Since $\kappa g$ too is continuous the latter term has the desired
continuity properties.

\underline{Remark.} From this proof it is clear that this is a result for 
\underline{harmonic} functions, rather that for holomorphic ones.

\section*{\protect\bigskip 3. Return to $H^{\infty }\left( \Sigma \right) $.
The space $H\left( \overline{X}\right) $.}

After this long digression we return back to the space $H^{\infty }\left(
\Sigma \right) $.$\medskip $

Let thus again $f\in H^{\infty }\left( \Sigma \right) $. From the discussion
in sec. 2 (after first rotating an angle $90^{\circ }$ in the opposite
sense) we conclude that the boundary values $f\left( i.\right) $ and $%
f\left( 1+i.\right) $ are well-defined distributions in the variable $y$
with values in $\Sigma $.$\medskip $

\underline{Remark.} For future use (see sec. 4) we remark that the functions
in $H^{\infty }\left( \Sigma \right) $ also take the \underline{same}
boundary values in the sense of the weak topology of $\Sigma $. Indeed let $%
f\in H^{\infty }\left( \Sigma \right) $. For any continuous linear
functional $l$ on $\Sigma $ ($l\in \Sigma ^{^{\prime }})$ the scalar
function $h\left( z\right) =l\left( f\left( z\right) \right) $ is in
(scalar) $H^{\infty }$. Thus $h$ takes (scalar) distributional boundary
values (and the latter agree of course with the classical pointwise a.e.
boundary values, known to exist by Fatou's theorem). It is clear that the
distributions $l\left( f\left( k+i.\right) \right) $ and $h\left(
k+i.\right) $ coincide.$\medskip $

Now we define $H^{\infty }\left( \overline{X}\right) $ to be the subspace of 
$H^{\infty }\left( \Sigma \right) $ of those functions $f\in H^{\infty
}\left( \Sigma \right) $ such that $f\left( i.\right) \in L^{\infty }\left(
X_{0}\right) $ and $f\left( 1+i.\right) \in L^{\infty }\left( X_{1}\right) .$
($L^{\infty }\left( X\right) $, $X$ any Banach space, is realized as a space
of vector valued distributions in the obvious way.\footnote{%
See note $\left\langle 3\right\rangle .$}) $H^{\infty }\left( \overline{X}%
\right) $ is a Banach space in the norm $\left\Vert f\right\Vert =\max
\left( \left\Vert f\left( i.\right) \right\Vert _{L^{\infty }\left(
X_{0}\right) },\left\Vert f\left( 1+i.\right) \right\Vert _{L^{\infty
}\left( X_{1}\right) }\right) $, indeed a Banach subspace of $H^{\infty
}\left( \Sigma \right) $. Both of these statements are embodied in the
inequality%
\begin{equation}
\left\Vert f\right\Vert _{H^{\infty }\left( \Sigma \right) }\leq \left\Vert
f\right\Vert _{H^{\infty }\left( \overline{X}\right) }.  \tag{1}
\label{eq-3-1}
\end{equation}%
Let us prove $\left( \ref{eq-3-1}\right) $. Since $L^{\infty }\left(
X_{k}\right) \subseteq L^{\infty }\left( \Sigma \right) $ $\left(
k=0,1\right) $ it is obvious that 
\begin{equation*}
\left\Vert f\right\Vert _{H^{\infty }\left( \overline{X}\right) }\geq \max
\left( \left\Vert f\left( i.\right) \right\Vert _{L^{\infty }\left( \Sigma
\right) },\left\Vert f\left( 1+i.\right) \right\Vert _{L^{\infty }\left(
\Sigma \right) }\right) .
\end{equation*}%
It thus suffices to prove that%
\begin{equation}
\left\Vert f\right\Vert _{H^{\infty }\left( \Sigma \right) }\leq \max
\left\{ \left\Vert f\left( i.\right) \right\Vert _{L^{\infty }\left( \Sigma
\right) },\left\Vert f\left( 1+i.\right) \right\Vert _{L^{\infty }\left(
\Sigma \right) }\right\}   \tag{$1^{^{\prime }}$}  \label{eq-3-1-prim}
\end{equation}%
that is, the maximum principle. The simplest way of establishing $\left( 
\text{\ref{eq-3-1-prim}}\right) $ is perhaps via the scalar case, in which
case we take $\left( \text{\ref{eq-3-1-prim}}\right) $ for granted. Let $l$
be any (continuous) linear functional on $\Sigma $ of norm $\leq 1$. Then
the scalar function $h\left( z\right) =l\left( f\left( z\right) \right) $,
as we have already noticed ultra (see the remark), certainly belongs to $%
H^{\infty }$ and in addition holds $\left\Vert h\left( k+i.\right)
\right\Vert _{L^{\infty }}\leq \left\Vert f\left( k+i.\right) \right\Vert
_{L^{\infty }\left( \Sigma \right) }$ $\left( k=0,1\right) $, since $%
\left\Vert l\right\Vert \leq 1$. It follows that $\left\Vert h\right\Vert
_{H^{\infty }}\leq $ the right hand side of $\left( \text{\ref{eq-3-1-prim}}%
\right) $. Since for any $z\in S$ holds $\left\Vert f\left( z\right)
\right\Vert _{\Sigma }=\sup \left\vert l\left( f\left( z\right) \right)
\right\vert $ $\left( \text{\ref{eq-3-1-prim}}\right) $ now follows.$%
\medskip $

Let further $A\left( \Sigma \right) $ be the space of those functions in $%
H^{\infty }\left( \Sigma \right) $ which take continuous boundary values
(belonging to $\Sigma )$. By the lemma in Sec. 2 this is to say that $f$ has
a continuous extension to $S$. We can then also define $A\left( \overline{X}%
\right) $ to be the subspace of $A\left( \Sigma \right) $ such that $f\left(
i.\right) $ is a continuous function with values in $X_{0}$ and $f\left(
1+i.\right) $ a continuous function with values in $X_{1}$, in symbols: $%
f\left( i.\right) \in C\left( X_{0}\right) $,$f\left( 1+i.\right) \in
C\left( X_{1}\right) $. Obviously $A\left( \overline{X}\right) $ is a closed
subspace of $H^{\infty }\left( \overline{X}\right) $.

\section*{4. A new characterization of the Calderon(-Lions) space $\overline{%
X}_{\protect\theta }$.}

By definition the image in $\Sigma $ of $A\left( \overline{X}\right) $ under
the evaluation map $f\mapsto f\left( \theta \right) $ is the (first)
Calderon space $\overline{X}_{\theta }$. (Here $\theta \in \left( 0,1\right) 
$; we identify of course $\left( 0,1\right) $ with a closed subset of $S^{i}$%
.$)$ As is well-known $\overline{X}_{\theta }$ is a Banach space in the
natural quotient norm.$\medskip $

Now we verify that we get the same space $\overline{X}_{\theta }$ if we in
this definition substitute $H^{\infty }$ for $A$.

\underline{Theorem.} The image of the evaluation map $H\left( \overline{X}%
\right) \rightarrow \Sigma :f\mapsto f\left( \theta \right) $ $\left( \theta
\in \left( 0,1\right) \right) $ is $\overline{X}_{\theta }$.

\underline{Proof: }. For convince, let us denote provisionally the image of $%
H\left( \overline{X}\right) $ by $\overline{X}_{\theta }^{\infty }$. Thus
our concern is to show that $\overline{X}_{\theta }=\overline{X}_{\theta
}^{\infty }$, indeed with equality of norms, if we again agree to equip $%
\overline{X}_{\theta }^{\infty }$ with the quotient norm.

One direction is easy. Indeed if $a\in \overline{X}_{\theta }$ then by
definition $a=f\left( \theta \right) $ for some function $f\in A\left( 
\overline{X}\right) $. Since $A\left( \overline{X}\right) \subseteq
H^{\infty }\left( \overline{X}\right) $ this shows that $\overline{X}%
_{\theta }\subseteq \overline{X}_{\theta }^{\infty }$.

For the opposite direction we first recall Calderon's famous inequality: If $%
f\in A\left( \overline{X}\right) $ then%
\begin{equation}
\left\Vert f\left( \theta \right) \right\Vert _{\overline{X}_{\theta }}\leq
\exp \left( \int \log \left\Vert f\left( i.\right) \right\Vert _{X_{0}}d\pi
_{0}^{\theta }+\int \log \left\Vert f\left( 1+i.\right) \right\Vert
_{X_{1}}d\pi _{1}^{\theta }\right)   \tag{$1$}  \label{eq-4-1}
\end{equation}%
where $\pi _{k}^{\theta }$ stands for the \underline{harmonic measure} of $%
\delta _{K}S$ at the point $\theta $; $d\pi _{k}^{\theta }/dy$ is thus the
corresponding part of the \underline{Poisson kernel} (evaluated at $\theta
):d\pi _{k}^{\theta }/dy=P\left( \theta ,k+iy\right) $; explicit expressions
for the latter can be found in \cite{5-Cal} or in \cite{3-BeL0}, chap. 4.

Let now $a\in \overline{X}_{\theta }^{\infty }$ so that $a=f\left( \theta
\right) $ for some $f\in H^{\infty }\left( \overline{X}\right) $.

Set $f_{n}\left( z\right) =n\tint\nolimits_{0}^{1/n}f\left( z+it\right)
dt,n=1,2,..$. We claim $f_{n}\in A\left( \overline{X}\right) $. Indeed it is
at once clear that $f_{n}$ is at least a function in $H^{\infty }\left(
\Sigma \right) $, whose boundary values by our principle of permanence (sec.
2) are given by an analogous formula: $f_{n}\left( k+iy\right)
=n\tint\nolimits_{0}^{1/n}f\left( k+iy+it\right) dt$, $\left( k=0,1\right) $%
. But the latter fact shows that $f_{n}\left( k+iy\right) $ is a continuous
function with values in $X_{k}$: in symbols $f_{n}\left( k+i.\right) \in
C\left( X_{k}\right) $ $\left( k=0,1\right) .$ (In fact $f_{n}\left(
k+i.\right) $ is even in $Lip\left( X_{k}\right) $ with the Lipschitz
constant bounded by $O\left( n\right) $. Because we can write%
\begin{equation*}
f_{n}\left( k+i\left( y+h\right) \right) -f_{n}\left( k+iy\right) =n\left(
\int_{1/n}^{1/n+h}-\int_{0}^{h}f\left( k+i\left( y+t\right) \right)
dt\right) 
\end{equation*}%
for $h$ small). Thus by the lemma in sec. 2 $f$ must be in $A\left( \Sigma
\right) $ and so in $A\left( \overline{X}\right) $ too. Our claim is
substantiated. We infer now that $a_{n}=f_{n}\left( \theta \right) \in 
\overline{X}_{\theta }$. Next we apply $\left( \text{\ref{eq-4-1}}\right) $
to the difference $f_{n}-f_{m}$. Then by Fatou's lemma we see that $\left\{
a_{n}\right\} $ is a Cauchy sequence in $\overline{X}_{\theta }$. But $%
\overline{X}_{\theta }$ is complete and continuously embedded in $\Sigma $.
Since the limit of $\left\{ a_{n}\right\} $ in $\Sigma $ clearly is $%
a=f\left( \theta \right) $ we therefore finally get $a\in \overline{X}%
_{\theta }$. $\square \medskip $

As a first application of this result (cf. infra sec. 6) we get that we
still have the same space $\overline{X}_{\theta }$ if we assume that the
boundary values are taken in the sense of the weak topology.$\medskip $

\underline{Corollary} (generalizing partial results by Janson \cite{23-Jan},
see e.g. th. 27). Let $a\in \Sigma $ and assume $a=f\left( \theta \right) $
for some $f\in H^{\infty }\left( \Sigma \right) $ such that for every $l\in
\Sigma ^{^{\prime }}$ holds $\lim_{z\rightarrow k+iy}l\left( f\left(
z\right) \right) =l\left( g_{k}\left( y\right) \right) $ a.e. where $g_{k}$
is some function in $L^{\infty }\left( X_{k}\right) $ $\left( k=0,1\right) $%
. Then $a\in \overline{X}_{\theta }$.

\underline{Proof:} It suffices to apply the remark in sec. 3. $\square $

\underline{Remark.} In this connection recall that a function with values in
a Banach space is holomorphic if and only if it is "weakly holomorphic"
(Dunford's theorem; see e.g. \cite{20-HiPh}, p.93).

\section*{5. The second Calderon space $\overline{X}^{\protect\theta }$.}

We now turn to the analogous question for the second Calderon space $%
\overline{X}^{\theta }$.$\medskip $

First recall the definition: $a\in \overline{X}^{\theta }$ if and only if $%
a\in \Sigma $ and there exists a function $g$ in $B\left( \overline{X}%
\right) $ such that $g^{^{\prime }}\left( \theta \right) =a$. Again $g\in
B\left( \overline{X}\right) $ means that $g$ is holomorphic in $S^{i}$ and
continuous in $S$ taking values in $\Sigma $, with $\left\Vert g\left(
z\right) \right\Vert _{\Sigma }=O\left( y\right) $, the boundary values $%
g\left( k+i.\right) $ being continuous functions with values in $X_{k}$ (at
least mod $\Sigma )$ subject to the condition 
\begin{equation*}
\left\Vert g\left( k+iy+it\right) -g\left( k+iy\right) \right\Vert
_{X_{k}}\leq c\left\vert t\right\vert ,
\end{equation*}%
$t$ real, $k=0,1$.

We also denote by $H_{\infty }\left( \overline{X}\right) $ the space of
functions $f$ which we obtain if we in the definition of $H^{\infty }\left( 
\overline{X}\right) $ (sec. 3) replace the functor $L^{1}\left( -\right) $
by the functor $\hom \left( L^{1},-\right) $; if $X$ is any Banach space $%
\hom \left( L^{1},X\right) $ consist of all continuous linear maps from $%
L^{1}$ into $X$. (Clearly $L^{\infty }\left( X\right) \subseteq \hom \left(
L^{1},X\right) $ for any $X$ but the converse holds true if and only if $X$
has the RN (Radon-Nikodym) property; this is practically a definition; see 
\cite{14-Fav}\footnote{%
See note $\left\langle 4\right\rangle .$}.) Thus $f\in H_{\infty }\left( 
\overline{X}\right) $ means that $f$ is in $H^{\infty }\left( \Sigma \right) 
$ with $f\left( k+i.\right) \in \hom \left( L^{1},X_{k}\right) $ $\left(
k=0,1\right) $.$\medskip $

\underline{Theorem.} The image of the evaluation map $H_{\infty }\left( 
\overline{X}\right) \rightarrow \Sigma :f\mapsto f\left( \theta \right) $ is
the Calderon space $\overline{X}^{\theta }$.

\underline{Proof:} In analogy with the corresponding proof in sec. 4 we
denote the image in question by $\overline{X}_{\infty }^{\theta }$. We want
to show that $\overline{X}^{\theta }=\overline{X}_{\infty }^{\theta }$.

Let $a\in \overline{X}_{\infty }^{\theta }$ so that $a=f\left( \theta
\right) $ for some $f\in H_{\infty }\left( \overline{X}\right) $. Put $%
g\left( z\right) \overset{def}{=}\tint\nolimits_{z_{0}}^{z}f\left( \zeta
\right) d\zeta $, $z_{0}$ some fixed point in $S$. Because of the
analyticity the integral is independent of contour of integration, to be
taken inside $S$ of course.

\underline{Claim.} $g\in B\left( \overline{X}\right) $. (Then will follows
that $a\in \overline{X}^{\theta }$ since it is clear that $g^{^{\prime
}}\left( \theta \right) =a$.$)$

\underline{Proof} (of the claim): It is clear that $g$ is holomorphic in $%
S^{i}$ with values in $\Sigma $ and that $\left\Vert g\left( z\right)
\right\Vert _{\Sigma }=O\left( y\right) $. There remains to investigate the
boundary values of $g$. To fix the ideas consider those on $\delta _{0}S$ $%
\left( x=0\right) $ taking for convenience $z_{0}=0$. Since $g^{^{\prime }}=f
$ the permanence principle (see sec. 2) shows that the distributional
derivative of $g\left( iy\right) $ is $if\left( iy\right) $. That is, for
every test function $\phi $ holds%
\begin{equation}
\int g\left( iy\right) \phi ^{^{\prime }}\left( y\right) dy=-\int f\left(
iy\right) \phi \left( y\right) dy.  \tag{1}  \label{eq-5-1}
\end{equation}%
Now put%
\begin{equation*}
F\left( y\right) \overset{def}{=}\int_{0}^{y}f\left( i\eta \right) d\eta
=\int_{-\infty }^{\infty }\chi _{\left( 0,y\right) }\left( \eta \right)
f\left( i\eta \right) d\eta .
\end{equation*}%
(If $\left( a,b\right) ,a<b$ is an interval on the real line we let $\chi
_{\left( a,b\right) }$ be its characteristic function; if $a>b$ we set $\chi
_{\left( a,b\right) }=-\chi _{\left( b.a\right) }$.) Then $F\left( 0\right)
=0$ and%
\begin{equation}
F\left( y+t\right) -F\left( y\right) =\int_{-\infty }^{\infty }\chi _{\left(
y,y+t\right) }\left( \eta \right) f\left( i\eta \right) d\eta .  \tag{2}
\label{eq-5-2}
\end{equation}%
Since by assumption $f\left( i.\right) \in \hom \left( L^{1},X_{0}\right) $
we see that $F\in C\left( X_{0}\right) $. Moreover by our choice of $z_{0}$
it is easy to see that \thinspace $g\left( iy\right) =-F\left( t\right) $.
(Take $\phi \left( y\right) =\tint\nolimits_{-\infty }^{y}\psi \left( \eta
\right) d\eta $ in $\left( \text{\ref{eq-5-1}}\right) $, where $%
\tint\nolimits_{-\infty }^{\infty }\psi \left( \eta \right) d\eta =0$.$)$ So 
$g\left( i.\right) $ is a continuous function with values in $X_{0}$. Using $%
\left( \text{\ref{eq-5-2}}\right) $ and once more $f\left( i.\right) \in
\hom \left( L^{1},X_{0}\right) $ we get the crucial estimate 
\begin{equation*}
\left\Vert g\left( iy+it\right) -g\left( iy\right) \right\Vert _{X_{0}}\leq
C\left\Vert \chi _{\left( y,y+t\right) }\right\Vert _{L^{1}}\leq C\left\vert
t\right\vert .
\end{equation*}%
In the same way we treat the boundary values on $\delta _{1}S$ $\left(
x=1\right) $. Therefore we can safely contend that $g\in B\left( \overline{X}%
\right) $. $\square $

\underline{Proof} (of the theorem/continued/): We have thus shown that $%
\overline{X}_{\infty }^{\theta }\subset \overline{X}^{\theta }$. To prove
the opposite inclusion consider any $a\in \overline{X}^{\theta }$ and let $g$
be in $B\left( \overline{X}\right) $ such that $g^{^{\prime }}\left( \theta
\right) =a$. Set $f=g^{^{\prime }}$.

\underline{Claim.} $f\in H_{\infty }\left( \overline{X}\right) $. (From
which at once follows that $a\in \overline{X}_{\infty }^{\theta }$ since $%
f\left( \theta \right) =g^{^{\prime }}\left( \theta \right) =a$.$)$

\underline{Proof} (of the claim): It is clear that $f$ is a holomorphic
function with values in $\Sigma $. Since $g\in B\left( \overline{X}\right) $
there holds in particular $\left\Vert g\left( k+iy+it\right) -g\left(
k+iy\right) \right\Vert _{\Sigma }\leq C\left\vert t\right\vert $ $\left(
k=0,1\right) $. Therefore by the maximum principle $\left\Vert g\left(
z+it\right) -g\left( z\right) \right\Vert _{\Sigma }\leq C\left\vert
t\right\vert $ $\left( z\in S\right) $. $f\left( z\right) =g^{^{\prime
}}\left( z\right) =\lim_{t\rightarrow 0}\left( g\left( z+it\right) -g\left(
z\right) \right) /it$ in the norm topology of $\Sigma $. Therefore $%
\left\Vert f\left( z\right) \right\Vert _{\Sigma }\leq C$ and $f\in
H^{\infty }\left( \Sigma \right) $.

Let us investigate the boundary values on $\delta _{0}S$, say. Again by the
permanence principle the distributional derivative of $g\left( iy\right) $
is $if\left( iy\right) $. That is, $\left( \ref{eq-5-1}\right) $ holds as
before.

But the left hand side is the limit in $X_{0}$ of 
\begin{eqnarray*}
&&\int_{-\infty }^{\infty }g\left( iy\right) \frac{\phi \left( y+t\right)
-\phi \left( y\right) }{t}dy \\
&=&\int_{-\infty }^{\infty }\frac{g\left( i\left( y-t\right) \right)
-g\left( iy\right) }{t}\phi \left( y\right) dy.
\end{eqnarray*}%
$g\in B\left( \overline{X}\right) $ again shows that the norm (in $X_{0}$)
of the latter expression for any $t$ is majorized by $C\left\Vert \phi
\right\Vert _{L^{1}}$. Therefore $\left\Vert \tint f\left( iy\right) \phi
\left( y\right) dy\right\Vert _{X_{0}}\leq C\left\Vert \phi \right\Vert
_{L^{1}}$ and $f\left( i.\right) \in \hom \left( L^{1},X_{0}\right) $.

In exactly the same way we handle the boundary values on $\delta _{1}S$ so
we get $f\left( 1+i.\right) \in \hom \left( L^{1},X_{1}\right) $ too. This
proves $f\in H_{\infty }\left( \overline{X}\right) $. $\square $

The proof is now complete. $\square $

\section*{6. Applications to the Calderon theory.}

\underline{N.B.} - The reading of this sec. is not needed for sec. 7-9.$%
\medskip $

\underline{6.1.} Calderon \cite{5-Cal} proved that if one of the spaces, say 
$X_{0}$, is reflexive then with equality of norms $\overline{X}_{\theta }=%
\overline{X}^{\theta }$ for any $\theta \in \left( 0,1\right) $. We now
generalize this result.

\underline{Proposition.} Assume that $X_{0}$ has the RN property, that is $%
\hom \left( L^{1},X_{0}\right) =L^{\infty }\left( X_{0}\right) $ if the
natural identifications are made. Then $\overline{X}_{\theta }=\overline{X}%
^{\theta }$ with equality of norms.

\underline{Remark.} This is generalization because by a theorem of
Phillips's (see \cite{11-DiUh}, p. 82) reflexive spaces have the RN
property. Notice also by a theorem of Bergh's \cite{2-BE79} we know that $%
\overline{X}_{\theta }$ always is \underline{isometrically} a subspace of $%
\overline{X}^{\theta }$.

If $X_{1}$ to has the RN property then $H^{\infty }\left( \overline{X}%
\right) =H_{\infty }\left( \overline{X}\right) $ and there is not much to
prove. For the general case we need again Calderon's inequality (see formula %
\ref{eq-4-1} in sec. 4).

\underline{Proof} (of the proposition): Let $a\in \overline{X}^{\theta }$
and pick up $f\in H_{\infty }\left( \overline{X}\right) $ such that $f\left(
\theta \right) =a$. (Here we use the theorem in sec. 5.) As in a previous
proof (sec. 4, theorem) we set $f_{n}\left( z\right)
=n\tint\nolimits_{0}^{1/n}f\left( z+it\right) dt$ $\left( n=1,2,...\right) $
or more properly $f_{n}\left( z\right) =n\tint \chi _{(y,y+1/n}\left( \eta
\right) f\left( x+i\eta \right) d\eta $ where we as usual write $z=x+iy$. It
is easily seen that $f\in H^{\infty }\left( \overline{X}\right) $ and even $%
\in A\left( \overline{X}\right) $. Therefore $a_{n}\overset{def}{=}%
f_{n}\left( \theta \right) \in \overline{X}_{\theta }$.

Because of the $RN$ property of the space $X_{0}$ we have $f\left( i.\right)
\in L^{\infty }\left( X_{0}\right) $ and $f_{n}\left( iy\right) \rightarrow
f\left( iy\right) $ a.e. as $n\rightarrow \infty $. Also $\left\Vert
f_{n}\left( k+i.\right) \right\Vert _{L^{\infty }\left( X_{k}\right) }\leq C$
$\left( k=0,1\right) $ for some $C$. Therefore using now Calderon's
inequality (formula \ref{eq-4-1} of sec. 4) and Fatou's lemma, in a manner
already familiar to us, we conclude that $\left\{ a_{n}\right\} $ is a
Cauchy sequence in $\overline{X}_{\theta }$ and from this again $a\in 
\overline{X}_{\theta }$. This proves $\overline{X}^{\theta }\subseteq 
\overline{X}_{\theta }$.The reverse inclusion is trivial.$\ \square \medskip 
$

\underline{Remark.} \underline{On the functor $\hom \left( L^{1},-\right) $}%
. Let $f\in \hom \left( L^{1},X\right) $ where $X$ is any given Banach
space; that is, $f$ is a continuous linear mapping from $L^{1}$ into $X$; $%
L^{1}$ (and $L^{\infty })$ is now taken with respect to some fixed measure
space $\Omega $ with measure $\mu $. The adjoint operator $f^{t}$ maps $%
X^{^{\prime }}$into $L^{\infty }$. For every element $l\in X^{^{\prime }}$%
with $\left\Vert l\right\Vert \leq 1$ we consider the element $f^{t}\left(
l\right) =l\circ f$ of $L^{\infty }$. Let now $f_{X}^{\#}$ denote the
supremum of all the $f^{t}\left( l\right) $ considered as element of the
Riesz space of (Radon) measure on $\Omega $ (cf. e.g. \cite{4-Bou}). This
supremum clearly exists and it is easily seen that it is a measure
absolutely continuous with respect to the given measure $\mu $. So $%
f_{X}^{\#}$ can be identified with an element of $L^{\infty }$ and one
readily verifies that $\left\Vert f\right\Vert =\left\Vert
f_{X}^{\#}\right\Vert _{L^{\infty }}$. In order to demonstrate the
usefulness of this notion let us indicate a quick proof of one of the
simplest positive results of the RN property, namely Dunford's theorem to
the effect that every Banach space with a \underline{conditionally bounded}
basis $\left( u_{n}\right) _{1}^{\infty }$ has the RN property (see \cite%
{11-DiUh}, p. 64). With no loss of generality we can assume that $\left(
u_{n}\right) _{1}^{\infty }$ in addition is \underline{monotone}. Let $%
\left( v_{n}\right) _{1}^{\infty }$be the dual "basis" in $X^{^{\prime }}$.
For $\phi \in L^{1}$ we can write $f\left( \phi \right) =\dsum_{1}^{\infty
}v_{n}\left( f\left( \phi \right) \right) u_{n}=\dsum_{1}^{\infty
}f^{t}\left( v_{n}\left( \phi \right) \right) u_{n}$. Each $f^{n}v_{n}$ can
be identified with a function in $L^{\infty }$. Consider $%
\dsum_{1}^{N}v_{n}\left( f\left( \phi \right) \right)
u_{n}=\tint\nolimits_{\Omega }\dsum_{1}^{N}f^{t}v_{n}\left( \omega \right)
u_{n}\phi \left( \omega \right) d\mu \left( \omega \right)
=\tint\nolimits_{\Omega }f_{N}\left( \omega \right) \phi \left( \omega
\right) d\mu \left( \omega \right) $. By the monotonicity we see that $%
\left\Vert f_{N}\left( \omega \right) \right\Vert _{X}\leq f_{X}^{\#}\left(
\omega \right) $ a.e. Therefore by the conditional boundedness $f_{N}\left(
\omega \right) \rightarrow f\left( \omega \right) $ a.e. for some function $%
f\in L^{\infty }\left( X\right) $. Clearly $f$ "represents" $f$ and the
proof is complete. \underline{Perhaps there are other applications too} (to
the RN theory). Here we content to point out returning to the arena of
complex interpolation that with the aid of this notion we can formally
generalize Calderon's inequality to the space $\overline{X}^{\theta }$ too.
If $f\in H_{\infty }\left( \overline{X}\right) $ then holds%
\begin{equation*}
\left\Vert f\left( \theta \right) \right\Vert _{\overline{X}^{\theta }}\leq
\exp \left( \int \log f_{X_{0}}^{\#}(i.\right) d\pi _{0}^{\theta }+\int \log
f_{X_{1}}^{\#}\left( 1+i.\right) d\pi _{1}^{\theta }
\end{equation*}%
with the same meaning of $\pi _{K}^{\theta }$ as in sec. 4 $\left(
k=0,1\right) $;  $f_{X_{k}}^{\#}\left( k+i.\right) $ is of course formed
with respect to Lebesque measure considering $f\left( k+i.\right) $ as an $%
X_{k}$ valued function.

The proof is about the same as for the classical Calderon's inequality in
the case of the space $\overline{X}_{\theta }$ (see \cite{5-Cal}, p. 134);
instead of $\left\Vert f\left( k+iy\right) \right\Vert _{X_{k}}$, which does
not "exists", use $f_{X_{k}}^{\#}\left( k+iy\right) $. Unfortunately we know
of no application of this new inequality.$\medskip $

\underline{6.2}. Having returned to the Calderon theory for good, we sketch
a new and \underline{conceptually} perhaps simpler proof of the duality
theorem (see \cite{5-Cal} or \cite{3-BeL0}), to the effect that $\left( 
\overline{X}_{\theta }\right) ^{^{\prime }}\cong \left( \overline{%
X^{^{\prime }}}\right) ^{\theta }$ if the dual pair $\overline{X^{^{\prime }}%
}=\left( X_{0}^{^{\prime }},X_{1}^{^{\prime }}\right) $ "exists". It will be
based on the following representation for the dual of $L^{1}\left( X\right) $%
, $X$ any Banach space: $\left( L^{1}\left( X\right) \right) ^{^{\prime
}}\cong \hom \left( L^{1},X^{^{\prime }}\right) $, which is different from
the one used by Calderon (and, what is important, not tied to the properties
of the real line).$\medskip $

Let thus $l$ be a continuous linear functional on $\overline{X}_{\theta }$.
We wish to identify $l$ with an element of $\left( \overline{X}^{^{\prime
}}\right) ^{\theta }$. (The other inclusion is trivial and will not be
considered here.) By Calderon's inequality once more the relation $f\mapsto
l\left( f\left( \theta \right) \right) $ defines a continuous linear
functional on $H^{\infty }\left( \overline{X}\right) $, \underline{equipped
with the norm} $\left\Vert \left\vert f\right\vert \right\Vert =\tint
\left\Vert f\left( iy\right) \right\Vert _{X_{0}}P\left( \theta ,iy\right)
dy+\tint \left\Vert f\left( 1+iy\right) \right\Vert _{X_{1}}P\left( \theta
,1+iy\right) dy$. So by the Hahn-Banach theorem and by the above
representation of the dual of $L^{1}\left( X\right) $ we see that there
exists elements $h_{0}$ and $h_{1}$ of $\hom \left( L^{1},X_{0}^{^{\prime
}}\right) $ and $\hom \left( L^{1},X_{1}^{^{\prime }}\right) $ respectively
such that (formally)%
\begin{equation*}
l\left( f\left( \theta \right) \right) =\int \left\langle h_{0}\left(
y\right) ,f\left( iy\right) \right\rangle P\left( \theta ,iy\right) dy+\int
\left\langle h_{1}\left( y\right) ,f\left( 1+iy\right) \right\rangle P\left(
\theta ,1+iy\right) dy.
\end{equation*}%
The proof is completed by observing that $h_{0}$ and~$h_{1}$ are the
(distributional) boundary values on $\delta _{0}S$ and $\delta _{1}S$
respectively of a suitable holomorphic function~$g$ with values in $\Sigma
\left( \overrightarrow{X^{^{\prime }}}\right) =X_{0}^{^{\prime
}}+X_{1}^{^{\prime }}$. This is done more or less as in Calderon's case (see 
\cite{5-Cal} or \cite{3-BeL0}) and one then also finds that $g\in H_{\infty
}\left( \overline{X^{^{\prime }}}\right) $.\ So we get $l\left( x\right)
=\left\langle y,x\right\rangle $ where $y=g\left( \theta \right) $ thus is
an element of $\left( \overline{X}^{^{\prime }}\right) ^{\theta }$.$\medskip 
$

\underline{6.3.} We next turn our attention to Janson's characterization (%
\cite{23-Jan}, th. 22) of the complex interpolation functors as orbit
functors in the sense of the Aronszajn-Gagliardo theorem \cite{1-AG65} (cf. 
\cite{3-BeL0}, chap. 2). Here we offer an alternative proof which again from
the \underline{conceptual point of view} might have some advantages. (What
we have in mind is of course possible extensions. E.g. the same procedure
should be applicable in the case of infinitely many spaces; cf.
Introduction.)$\medskip $

First we reformulate Janson's theorem in a way suitable for our purposes.
(Janson uses the discrete version of the Calderon spaces (Cwikel's theorem 
\cite{10-Cwi}) but this is not really the point.)$\medskip $

Consider the following (compatible) pair of Banach spaces $\overline{F}%
=\left( F_{0},F_{1}\right) \footnote{%
See note $\left\langle 5\right\rangle .$}$. The "containing" space is simply
the dual of (scalar) $H^{\infty }$ and the space $F_{k}$ consists of those
elements $\mu $ of $\left( H^{\infty }\right) ^{^{\prime }}$ which can be
represented in the form $\mu \left( \phi \right) =\tint\nolimits_{.\infty
}^{\infty }\psi \left( y\right) \phi \left( k+iy\right) dy$, $\left( \phi
\in H^{\infty }\right) $ with $\psi \in L^{1}$ $\left( k=0,1\right) $.$%
\medskip $

It follows that an element $\mu \in \left( H^{\infty }\right) ^{^{\prime }}$
is in the hull (sum) $F_{0}+F_{1}$ if and only if it can be represented in
the form $\mu \left( \phi \right) =\tint\nolimits_{-\infty }^{\infty }\psi
_{0}\left( y\right) \phi \left( iy\right) dy+\tint\nolimits_{.\infty
}^{\infty }\psi _{1}\left( y\right) \phi \left( 1+iy\right) dy$ with $\psi
_{k}\in L^{1}$ $\left( k=0,1\right) $. Let further $F=F_{\theta }$ $\left(
\theta \in \left( 0,1\right) \right) $ denote the subspace of $\left(
H^{\infty }\right) ^{^{\prime }}$ consisting of the linear functionals of
the form $\mu \left( \phi \right) =\tint\nolimits_{-\infty }^{\infty }\psi
\left( y\right) \phi \left( \theta +iy\right) dy$ with $\psi \in L^{1}$ and
by $\delta _{z}$ $\left( z\in S\right) $ the linear functional defined by $%
\delta _{z}\left( \phi \right) =\phi \left( z\right) $ (in other words, the
evaluation map).$\medskip $

Then Janson's results says that

$1^{\circ }$ $\overline{X}_{\theta }$ is the orbit of $F$ in $\overline{X}$.

$2^{\circ }$ $\overline{X}^{\theta }$ is the orbit of $\delta _{\theta }$ in 
$\overline{X}$.

($\overline{X}$ denotes as before an arbitrary Banach couple.).

\underline{Proof:} . One way is easy. By general principles connected with
the Aronszajn-Gagliardo theorem \cite{1-AG65} it suffices to show $%
F\subseteq \overline{F}_{\theta }$ and $\delta _{\theta }\in \overline{F}%
^{\theta }$ respectively.

In the first case if $\mu \in F$, $\mu \left( \phi \right) =\tint \psi
\left( y\right) \phi \left( \theta +iy\right) dy$ $\left( \phi \in H^{\infty
}\right) $ for some $\psi \in L^{1}$, we get a vector valued function $f$ by
defining $f\left( z\right) $ for $z\in S^{i}$ as the linear functional $\phi
\mapsto \tint \psi \left( \eta \right) \phi \left( z+i\eta \right) d\eta $.
It is readily seen $f\in H^{\infty }\left( \overline{F}\right) $ and that $%
f\left( \theta \right) =\mu $. This proves $\mu \in \overline{F}_{\theta }$.

In the second case set $f\left( z\right) =\delta _{z}$. Then $f\in H_{\infty
}\left( \overline{F}\right) $ and trivially $f\left( \theta \right) =\delta
_{\theta }$ which gives $\delta _{\theta }\in \overline{F}^{\theta }$. The
fact that $f\in H_{\infty }\left( \overline{F}\right) $ requires a proof.
Let us just indicate how one sees that $f\left( i.\right) \in \hom \left(
L^{1},F_{0}\right) $. The point is that the pointwise boundary values $%
\delta _{iy}$ do not belong to the sum $F_{0}+F_{1}$. But if we smear them
out with a test function $\psi $ in we get the linear functional $\phi
\mapsto \tint\nolimits_{-\infty }^{\infty }\psi \left( y\right) \phi \left(
iy\right) dy$ $\left( \phi \in H^{\infty }\right) $ which obviously is in $%
F_{0}$. This fixes the matter.

Now we turn to the opposite inclusion. If $a$ is $1^{\circ }$ in $\overline{X%
}_{\theta }$ or $2^{\circ }$ in $\overline{X}^{\theta }$ let $f$ be in $%
H^{\infty }\left( \overline{X}\right) $ or in $H_{\infty }\left( \overline{X}%
\right) $ respectively with $a=f\left( \theta \right) $. First define a
linear mapping $U:\overline{F}\rightarrow \overline{X}$ as follows. If $\mu
\in \Sigma \left( \overline{F}\right) =F_{0}+F_{1}$, $\mu \left( \phi
\right) =\tint\nolimits_{-\infty }^{\infty }\psi _{0}\left( y\right) \phi
\left( iy\right) dy+\tint\nolimits_{-\infty }^{\infty }\psi _{1}\left(
y\right) \phi \left( 1+iy\right) dy$ $\left( \phi \in H^{\infty }\right) $
with $\psi _{k}\in L^{1}$ $\left( k=0,1\right) $ we set $U\left( \mu \right)
=\tint\nolimits_{-\infty }^{\infty }\psi _{0}\left( y\right) f\left(
iy\right) dy+\tint\nolimits_{-\infty }^{\infty }\psi _{1}\left( y\right)
f\left( 1+iy\right) dy$. This definition makes sense, since we know at least
that $f\left( k+iy\right) \in \hom \left( L^{1},X_{k}\right) $ $\left(
k=0,1\right) $. It is also practically obvious that $U\left( \mu \right) $
does not depend on the particular representation of $\mu $ in terms of
functions $\psi _{0}$ and $\psi _{1}$. Moreover clearly $U:F_{k}\rightarrow
X_{k}$ $\left( k=0,1\right) $, that is, in abbrevia $U:\overline{F}%
\rightarrow \overline{X}$.

In case $2^{\circ }$ we simply put $\mu =\delta _{\theta }$. Then we can
take $\psi _{k}\left( y\right) =P\left( \theta ,k+iy\right) $ (Poisson
kernel) so $U\left( \delta _{\theta }\right) =f\left( \theta \right) =a$ and 
$a$ is in the orbit of $\delta _{\theta }$.

In the case $1^{\circ }$ we must (as in \cite{23-Jan}) use an approximation
device. Set $a_{n}=n\tint\nolimits_{0}^{1/n}f\left( \theta +it\right) dt$.
Then each $a_{n}\in \overline{X}_{\theta }$ and $a_{n}$ tends to $a$ in $%
\overline{X}_{\theta }$ (Here we used $f\left( k+iy\right) \in L^{\infty
}\left( X_{k}\right) $.$)$. Also $a_{n}=U\left( \mu _{n}\right) $ where $\mu
_{n}$ is the linear functional defined by $\mu _{n}\left( \phi \right)
=n\tint\nolimits_{0}^{1/n}\phi \left( \theta +it\right) dt$. Since, as is
readily seen, $\mu _{n}\in F$ each $a_{n}$ lies in the orbit of $F$ and its
(orbit) norm can be uniformly bounded by the one of $a$ (in $\overline{X}%
_{\theta })$. This shows (apply the usual iterative procedure) that $a$ too
is in the orbit.$\medskip $

In conclusion let us point out that as a consequence of his theorem Janson
also obtains a simple proof of the reiteration theorem (see \cite{23-Jan},
th. 25).

\section*{7. A vector valued (or better Banach bundle valued) interpolation
problem.}

Let $\left\{ z_{n}\right\} _{n=1}^{\infty }\subseteq S^{i}$ be any (fixed)
sequence \underline{uniformly separating} in the sense of Carleson \cite%
{6-Car}. Then by Carleson's theorem \cite{6-Car} for any (scalar) sequence $%
\left\{ w_{n}\right\} _{n=1}^{\infty }\subseteq l^{\infty }$ we can find a
(scalar) function $f\in H^{\infty }$ such that $f\left( z_{n}\right) =w_{n}$ 
$\left( n=1,2....\right) $.$\medskip $

Now we turn to the corresponding vector valued problem: Given a (vector
valued) sequence $\left\{ w_{n}\right\} _{n=1}^{\infty }$ to find a function 
$f\in H^{\infty }\left( \overline{X}\right) $ such that $f\left(
z_{n}\right) =w_{n}$ $\left( n=1,2,...\right) $. Clearly a \underline{%
necessary} condition for this to be possible is that $w_{n}\in \overline{X}%
_{z_{n}}$ $\left( n=1,2,...\right) $, $\sup_{n}\left\Vert w_{n}\right\Vert _{%
\overline{X}_{z_{n}}}<\infty $. ($\overline{X}_{z}$ $\left( z\in
S^{i}\right) $ denotes of course the space $\overline{X}_{Rez}$.$)\medskip $

\underline{Claim.} This necessary condition is \underline{sufficient} too.$%
\medskip $

To substantiate this claim we will make use of the functions constructed by
Per Beurling (see \cite{7-Car}): There exists a sequence of (scalar)
functions $\left\{ F_{n}\right\} _{n=1}^{\infty }\subseteq H^{\infty }$ such
that $F_{n}\left( z_{k}\right) =\delta _{nk}$ $\left( n,k=1,2,...\right) $
and $\dsum_{n=1}^{\infty }\left\vert F_{n}\left( z\right) \right\vert \leq
C<\infty $ $\left( z\in S^{i}\right) $. Then in the scalar case a particular
solution of our interpolation problem is provided by%
\begin{equation}
f\left( z\right) =\dsum_{n=1}^{\infty }w_{n}F_{n}\left( z\right) .  \tag{$1$}
\label{eq-7-1}
\end{equation}%
In the vector valued case we must modify $\left( \text{\ref{eq-7-1}}\right) $
somewhat. Indeed by definition (and by the theorem in sec. 4) we can find
functions $f_{n}\in H^{\infty }\left( \overline{X}\right) $ such that $%
f_{n}\left( z_{n}\right) =w_{n}$ $\left( n=1,2,...\right) $ and such that $%
C^{^{\prime }}\overset{def}{=}\sup_{n}\left\Vert f_{n}\right\Vert
_{H^{\infty }\left( \overline{X}\right) }<\infty $. Then we put 
\begin{equation}
f\left( z\right) =\dsum_{n=1}^{\infty }f_{n}\left( z\right) F_{n}\left(
z\right) .  \tag{$1^{^{\prime }}$}  \label{eq-7-1-prim}
\end{equation}

\underline{Sub-claim.} The function $f$ defined by $\left( \text{\ref%
{eq-7-1-prim}}\right) $ is in $H^{\infty }\left( \overline{X}\right) $ and
solves our interpolation problem, i.e. $f\left( z_{n}\right) =w_{n}$ $\left(
n=1,2,...\right) $.

\underline{Proof} (of the sub-claim). It is clear that for any $z\in S^{i}$
holds $\left\Vert f_{n}\left( z\right) \right\Vert _{\Sigma }\leq \left\Vert
f_{n}\right\Vert _{H^{\infty }\left( \Sigma \right) }\leq \left\Vert
f_{n}\right\Vert _{H^{\infty }\left( \overline{X}\right) }\leq C^{^{\prime }}
$ $\left( n=1,2,...\right) $. Therefore the series $\dsum_{n=1}^{\infty
}f_{n}\left( z\right) F_{n}\left( z\right) $ is normally convergent in $%
\Sigma $ for $z\in S^{i}$ and its sum $f\left( z\right) $ is clearly a
holomorphic function with $\left\Vert f\left( z\right) \right\Vert _{\Sigma
}\leq CC^{^{\prime }}$. In particular at least $f\in H^{\infty }\left(
\Sigma \right) $ and $f\left( z_{n}\right) =w_{n}$ $\left( n=1,2,...\right) $%
.$\medskip $

It remains as usual to investigate the boundary values. By the principle of
permanence (sec. 2) we have for the boundary functions on $\delta _{0}S$,
say, 
\begin{eqnarray*}
f\left( i.\right)  &=&\dsum f_{n}\left( i.\right) F_{n}\left( i.\right) , \\
\dsum \left\vert F_{n}\left( i.\right) \right\vert  &\leq &C<\infty
,\left\Vert f_{n}\left( i.\right) \right\Vert _{X_{0}}\leq C.
\end{eqnarray*}%
Therefore $f\left( i.\right) \in L^{\infty }\left( X_{0}\right) $. Similarly
we prove that $f\left( 1+i.\right) \in L^{\infty }\left( X_{1}\right) $. So
that indeed $f\in H^{\infty }\left( \overline{X}\right) $. $\square $

\section*{8. A Banach bundle valued $\overline{\protect\delta }$ problem.}

As is well-known (see \cite{21-Hor}) the interpolation problem in the
previous sec. is a special case of a problem for the $\overline{\delta }$%
-operator which in the scalar case can be formulated as follows: Let $\mu $
be a Carleson measure (see \cite{6-Car}, \cite{21-Hor}, \cite{24-Jon}), to
fix the ideas \underline{positive}, and consider any function $w\in
L^{\infty }\left( \mu \right) $. To find a function $u$ in $S^{i}$ such that 
$\overline{\delta }u=w\mu $ $\left( \overline{\delta }=\delta /\delta 
\overline{z}\right) $ and such that its distributional boundary values (in
the sense of sec. 2) lies in $L^{\infty }$. (To see the connection with the
interpolation problem of sec. 7 we remark that if $\mu $ is discrete then $u$
is related to the previous $f$ by $u=f/B$ where $B$ is the Blaschke product
formed with the sequence $\left\{ z_{n}\right\} _{n=1}^{\infty }$.$)$ We can
assume that at least $u\in L^{1}$, in which case we know for sure (sec. 2)
that the boundary distribution exists. $\medskip $

This problem always has a solution, in view of the Hahn-Banach theorem (see 
\cite{21-Hor}), but it is of course not unique. Once an appropriate
substitute for the Per Beurling functions has been found we can however
write down a particular solution. To this end let us consider a kernel $%
K\left( z,\zeta \right) $ assumed to be at least measurable in both
variables (measurable in $\zeta $ with respect to $\mu )$ such that

$\left( i\right) :K\left( z,\zeta \right) $ is analytic in $z$ for $\zeta $
fixed $\mu $ a.e.,

$\left( ii\right) :K\left( \zeta ,\zeta \right) =1$,

$\left( iii\right) :\tint K\left( z,\zeta \right) P\left( z,\zeta \right)
d\mu \left( \zeta \right) \leq C<\infty $ $\left( z\in \delta S\right) $.$%
\medskip $

Here $P\left( z,\zeta \right) $ $\left( z\in \delta S,\zeta \in S\right) $
stands for the Poisson kernel. We define the kernel $H\left( z,\zeta \right) 
$ $\left( z\in S,\zeta \in S\right) $ by stipulating that $1^{\circ }$ it
should for $\zeta $ fixed be analytic in $z$ except at $z=\zeta $ where it
should have a simple pole with the residue $1/2\pi i$ and that $2^{\circ }$ $%
\left\vert H\left( z,\zeta \right) \right\vert =P\left( z,\zeta \right) $ $%
\left( z\in \delta S,\zeta \in S\right) $. (See \cite{24-Jon} where a fairly
explicit construction of such a kernel $K\left( z,\zeta \right) $ is given.%
\footnote{%
See note $\left\langle 6\right\rangle .$}) Now the desired solution is
obtained as%
\begin{equation}
u\left( z\right) =\int K\left( z,\zeta \right) H\left( z,\zeta \right)
w\left( \zeta \right) d\mu \left( \zeta \right) \text{.}  \tag{$1$}
\label{eq-8-1}
\end{equation}%
$\medskip $

Next let us consider the generalization to the vector valued case. It seems
natural to ask for a solution $u$ of $\overline{\delta }u=w\mu $ such that $%
u\left( i.\right) \in L^{\infty }\left( X_{0}\right) $, $u\left( 1+i.\right)
\in L^{\infty }\left( X_{1}\right) $. What do we then have to require from $%
w?$ (A priori we assume only that $w$ is $\Sigma $ valued.) A minute's
reflection gives that we should assume $w\left( \zeta \right) \in \overline{X%
}_{\zeta }$ $\mu -$a.e., $\sup \left\Vert w\left( \zeta \right) \right\Vert
_{\overline{X}_{\zeta }}<\infty $ - to be exact, sup here means $\mu -$ess
sup - along with some measurability condition too. We postpone for a moment
the discussion of the exact form of the latter. Since $w\left( \zeta \right)
\in \overline{X}_{\zeta }$ we can for each $\zeta $ ($\mu -a.e.)$ find a
function $W\left( \cdot ,\zeta \right) \in H^{\infty }\left( \overline{X}%
\right) $ such that $w\left( \zeta \right) =W\left( \zeta ,\zeta \right) $.
We can also assume that $\left\Vert W\left( \cdot ,\zeta \right) \right\Vert
_{\overline{X}_{\zeta }}\leq C^{^{\prime }}<\infty $ $\left( \zeta \in
S\right) $. Then we set (formally)%
\begin{equation}
u\left( z\right) =\int_{S}K\left( z,\zeta \right) H\left( z,\zeta \right)
W\left( z,\zeta \right) d\mu \left( \zeta \right) \text{ }\left( z\in
S^{i}\right) \text{.}  \tag{$2$}  \label{eq-8-2}
\end{equation}%
Note that this formula $\left( \text{\ref{eq-8-2}}\right) $ generalizes both 
$\left( \text{\ref{eq-8-1}}\right) $ above and $\left( \text{\ref%
{eq-7-1-prim}}\right) $ of sec. 7, the discrete case. If everything works
out smoothly we expect $\left( \text{\ref{eq-8-2}}\right) $ to provide us
with a solution of $\overline{\delta }y=w\mu $ of the desired kind. Indeed
by the properties $\left( i\right) -\left( ii\right) $ of $K\left( z,\zeta
\right) $ and by property $1^{\circ }$ of $H\left( z,\zeta \right) $ we have
(formally) $\overline{\delta }K\left( \cdot ,\zeta \right) H\left( \cdot
,\zeta \right) W\left( \cdot ,\zeta \right) =w\left( \eta \right) \delta
\left( \cdot -\zeta \right) ,$ so $\overline{\delta }u\left( z\right) =\tint
\delta \left( z-\zeta \right) w\left( \zeta \right) d\mu \left( \zeta
\right) =w\left( z\right) \mu \left( z\right) $. Also if we can pass to the
limit in $\left( \text{\ref{eq-8-2}}\right) $ we expect that the same
identity holds for $z\in \delta S$. From this using property $\left(
iii\right) $ and $2^{\circ }$ we infer that%
\begin{eqnarray*}
\left\Vert u\left( z\right) \right\Vert _{X_{k}} &\leq &\int \left\vert
K\left( z,\zeta \right) \right\vert P\left( z,\zeta \right) \left\Vert
W\left( z,\zeta \right) \right\Vert _{X_{k}}d\mu \left( \zeta \right)  \\
&\leq &CC^{^{\prime }}\text{ }\left( k=0,1,z\in \delta S\right) 
\end{eqnarray*}%
so that we get the correct boundary behavior. $\medskip $

To make all this rigorous we have to make precise the assumptions on $w$.
Let us denote by $\overline{\overline{X}}$ the family of Banach spaces $%
\left( \overline{X}_{\zeta }\right) _{\zeta \in S^{i}}$. Let us also
introduce the norm%
\begin{equation*}
\left\Vert w\right\Vert =\mu -\text{ess sup}_{\zeta \in S^{i}}\left\Vert
w\left( \zeta \right) \right\Vert _{\overline{X}_{\zeta }}\text{.}
\end{equation*}

\underline{Definition.} We say that $w$ is in $L^{\infty }\left( \mu ,%
\overline{\overline{X}}\right) $ is $w$ can be approximated in the preceding
metric $\left\Vert w\right\Vert $ with a sequence of functions $w_{n}$ of
the form 
\begin{equation*}
w_{n}\left( \zeta \right) =\dsum_{v=1}^{N_{n}}\chi _{e_{v}^{n}}\left( \zeta
\right) f_{v}^{n}\left( \zeta \right) 
\end{equation*}%
where each $e_{v}^{n}$ is a $\mu $ measurable subset of $S^{i}$, $\chi
_{e_{v}^{n}}$ standing for its characteristic function, and the $f_{v}^{n}$
belong to $H^{\infty }\left( \overline{X}\right) $.$\medskip $

With this definition it is an easy matter to prove rigorously the following
result.

\underline{Proposition.} For any $w\in L^{\infty }\left( \mu ,\overline{%
\overline{X}}\right) $ there exists a function $u\in L_{loc}^{1}\left(
\Sigma \right) $ such that $\overline{\delta }u=w\mu $ and such that in the
sense of distributions $u\left( k+i.\right) \in L^{\infty }\left(
X_{k}\right) $ $\left( k=0,1\right) $.

\underline{Proof:} We assume first that $w$ itself is of the form $w\left(
\zeta \right) =\dsum_{v=1}^{N}\chi _{e_{v}}\left( \zeta \right) f_{v}\left(
\zeta \right) $ where the $e_{v}$ are $\mu $ measurable sets and $f_{v}\in
H^{\infty }\left( \overline{X}\right) $.Then we set 
\begin{eqnarray*}
u\left( z\right)  &=&\dsum_{v=1}^{N}f_{v}\left( z\right) \int_{e_{v}}K\left(
z,\zeta \right) H\left( z,\zeta \right) d\mu \left( \zeta \right) = \\
&&\dsum_{v=1}^{n}f_{v}\left( z\right) u_{v}\left( z\right) \text{.}
\end{eqnarray*}%
(This is formally $\left( \text{\ref{eq-8-2}}\right) $ with $W\left( z,\zeta
\right) =\dsum \chi _{e_{v}}\left( \zeta \right) f_{v}\left( z\right) !)$ By
Jones's results \cite{24-Jon} we certainly have for each $v$ $\overline{%
\delta }u_{v}=\chi _{e_{v}}\mu $. Therefore we get $\overline{\delta }%
u=\dsum f_{v}\overline{\delta }u_{v}=\dsum f_{v}\chi _{e_{v}}\mu =w\mu $
(since $\overline{\delta }f_{v}=0!)$. It is clear that $u$ is in $%
L_{loc}^{1}\left( \Sigma \right) $ and satisfies all the requirements of the
proposition. In particular we have%
\begin{equation*}
\left\Vert u\left( k+i.\right) \right\Vert _{L^{\infty }\left( X_{k}\right)
}\leq C_{k}\left\Vert w\right\Vert \text{ }\left( k=0,1\right) \text{.}
\end{equation*}

Moreover it is easy to see that we have for each $R$ an estimate of the type%
\begin{equation*}
\int\nolimits_{\left\vert y\right\vert \leq R}\left\Vert u\right\Vert
_{\Sigma }dxdy\leq C_{R}\left\Vert w\right\Vert \text{.}
\end{equation*}%
It is now obvious how to treat the general case too. If $\left\{
w_{n}\right\} _{n=1}^{\infty }$ is an approximating sequence for $w$ in the
sense of the above definition, $\left\Vert w-w_{n}\right\Vert \rightarrow 0$
as $n\rightarrow \infty $ and if $\left\{ u_{n}\right\} $ is the
corresponding sequence of solutions of $\overline{\delta }u_{n}=w_{n}$
obtained by the previous procedure then $u_{n}$ tends to a solution $u$ of $%
\overline{\delta }u=w\mu $ solving our problem.$\medskip $

\underline{Remark.} A final comment is in order. The definition of $%
L^{\infty }\left( \mu ,\overline{\overline{X}}\right) $ ultra (as well as
the notation) might seem quite \underline{ad hoc}, but it is not; in fact,
we insist, it is the most natural thing to be thought of. To see this we
first have to change slightly our point of view. Instead of considering the
family of Banach spaces $\left( \overline{X}_{\zeta }\right) _{\zeta \in
S^{i}}$ (until now denoted by $\overline{\overline{X}})$ we consider their 
\underline{set theoretic union}, using for this object the previous symbol $%
\overline{\overline{X}}$. Then $\overline{\overline{X}}$ might be considered
a vector bundle over $S^{i}$ (with no topology yet, however): $\medskip $

We have a natural projection $\pi :\overline{\overline{X}}\rightarrow S^{i}$
and each fiber $\pi ^{-1}\left( \zeta \right) =\overline{X}_{\zeta }$ is a
Banach space, thus \underline{a fortiori} a vector space. But more, $%
\overline{\overline{X}}$ is in fact in a technical sense (see Fell \cite%
{15-Fel}, \cite{16-Fel}) a \underline{Banach bundle}. In particular $%
\overline{\overline{X}}$ thus indeed carriers a natural topology by its own
right. By a classical procedure by Godement's \cite{18-God}, \cite{19-God}
the Banach bundle structure can quite generally (any vector bundle over a
(usually) locally compact space $\Omega )$ be defined by first specifying a
suitable family of "principal sections". In our case $\left( \Omega
=S^{i}\right) $ there is a canonical choice of the principal sections. These
are simply the sections of the type $S^{i}\rightarrow \overline{\overline{X}}%
:\zeta \mapsto f\left( \zeta \right) $,  $f$ a function in $H^{\infty
}\left( \overline{X}\right) $. For any Banach bundle $B$ and any positive
measure $\mu $ over its base space $\Omega $ we can construct a theory of $%
L^{p}$ spaces $L^{p}\left( \mu ,B\right) $ (apart from the works already
listed see the books by Dixmier \cite{13-Dix}, p. 186-194 and by Dinculeanu 
\cite{12-Din}, p. 413-414 where additional references can be found). If we
specialize we see that our $L^{\infty }\left( \mu ,\overline{\overline{X}}%
\right) $ ($\mu $ is now again a Carleson measure) is just a special case of
these general spaces. We reader can also convince himself of the fact that
in the scalar case our definitions reduces to (one of) the usual definitions
of $L^{\infty }\left( \mu \right) $. Also it is now finally plain why the
word "Banach bundle" appears in the title of this and the previous sec.

\section*{9. Other problems for $H^{\infty }$.}

After this initial success it is now natural to pause and to ask what other
traditional problems for $H^{\infty }$can be generalized to this new setting
of Banach bundles. In the first place what comes to ones mind is the Corona
Problem but we have not obtained any positive results in that direction so
presently we are bound to think that this might be a quite hopeless thing.
One can also ask if one could do something with $H^{p}$. If we understand
correctly a remark in \cite{24-Jon} it should be quite easy to prove vector
valued analogous of the interpolation theorem of Shapiro and Shields \cite%
{26-ShSh} but we have not tried to carry out the details.

\section*{\protect\underline{Notes.}}

$\left\langle 1\right\rangle $. Convention: If $V$ is any normed space,
fixed under the discussion, we denote its norm by $\left\Vert \cdot
\right\Vert $. If there are several spaces involved, in order to avoid
confusion, we use $V$ as subscript, thus writing $\left\Vert \cdot
\right\Vert _{V}$ for $\left\Vert \cdot \right\Vert $.

$\left\langle 2\right\rangle $.We are in particular interested in the
(rather trivial) special case $A=0$ (bounded functions) but also in
derivatives of bounded functions. If $f$ is bounded then $f^{^{\prime
}}\left( z\right) =O\left( y^{-1}\right) $ by an easy application of
Cauchy's theorem.

$\left\langle 3\right\rangle $. Unless otherwise specified by $L^{\infty }$
we mean $L^{\infty }$ with respect to Lebesque measure on the real line.
Similarly for $L^{1}$.

$\left\langle 4\right\rangle $. This work is entitled "Vector measures" but
is in large portions devoted to a study of Banach space having the $RN$
property.

$\left\langle 5\right\rangle $. $F$ stands for Fourier.

$\left\langle 6\right\rangle $. In \cite{24-Jon} there is considered of
course not the strip $S$ but the upper (Poincare) half-plane $\Pi =\left\{
y\geq 0\right\} $, the whole set-up being essentially invariant for
biholomorphic (conformal) transformations; in particular we have the
canonical map $z\rightarrow \exp \left( i\pi z\right) $ of $S^{i}$ onto $\Pi
^{i}$. In the case of the half-plane $\Pi $ one has%
\begin{eqnarray*}
H\left( z,\zeta \right)  &=&\frac{1}{2\pi i}\frac{\zeta -\overline{\zeta }}{%
\left( z-\zeta \right) \left( z-\overline{\zeta }\right) }\text{ }\left(
z\in \Pi ^{i}\right) \text{,} \\
P\left( x,\zeta \right)  &=&\frac{1}{\eta i}\frac{\func{Im}\zeta }{%
\left\vert x-\zeta \right\vert ^{2}}\text{ }\left( x=z\in \delta \Pi \right) 
\text{.}
\end{eqnarray*}

\label{_end}

\end{document}